\documentclass[12pt]{article}
\title{}
\author{}
\usepackage{amsmath, amsfonts}
\def\be{\begin{equation}}
\def\ee{\end{equation}}
\def\bm{\mathbf}

\def\R{\mathbb R}
\def\Z{\mathbb Z}
\def\qed{\hfill$\diamondsuit$}
\def\ed{\end{document}}

\newtheorem{prop}{Proposition}

\newtheorem{cor}{Corollary}

\def\eqalign#1{\null\,\vcenter{\openup\jot\ialign
              {\strut\hfil$\displaystyle{##}$&$\displaystyle{{}##}$
               \hfil\crcr#1\crcr}}\,}

\title{{\bf A Central Limit Theorem For \\ Linear Random Fields}}
\author{Atul Mallik and Michael Woodroofe \\ The University of Michigan}
\date{}

\begin{document}
\maketitle

\section{Introduction}

	Consider a two-dimensional, linear random field, say
$$
	X_{j,k} = \sum_{r\in\Z}\sum_{s\in \Z} a_{r,s}\xi_{j-r,k-s} = \sum_{r\in\Z}\sum_{s\in \Z} a_{j+r,k+s}\xi_{-r,-s},
$$
where $a_{r,s},\ r,s \in \Z$, are square summable, $\xi_{r,s},\ r,s \in \Z$, are i.i.d. with mean $0$ and unit variance, and $\Z$ denotes the integers.  It is convenient to regard the array $\bm{a} = (a_{r,s}:r,s \in \Z)$ as an element of $\ell^2(\Z^2)$.  Let $F$ denote the common distribution function of the $\xi_{r,s}$ and $(\Omega,{\cal A},P)$ the probability space on which they are defined. 
 If $\Gamma$ is a finite subset of $\Z^2$, let 
$$
	S = S(\bm{a},\Gamma) = \sum_{(j,k)\in\Gamma} X_{j,k}
$$
and
$$
	\sigma^2 = \sigma^2(\bm{a},\Gamma) = E(S^2),
$$
and suppose that $\sigma^2 > 0$. (Of course, $S$ depends on $\omega \in \Omega$ too, but this dependence is suppressed.  As indicated dependence on $\bm{a}$ and $\Gamma$ will only be displayed when needed for clarity.)  Then 
$$
	S = \sum_{r\in\Z}\sum_{s\in \Z} b_{r,s}\xi_{-r,-s},
$$
where
$$
	b_{r,s} = b_{r,s}(\bm{a},\Gamma) =  \sum_{(j,k)\in \Gamma} a_{j+r,k+s},
$$
and
$$
	\sigma^2 = \sum_{r\in\Z}\sum_{s\in \Z} b_{r,s}^2,
$$
assumed to be positive.  Let $\Phi$ denote the standard normal distribution and 
$$
	G(z) = G(z;\bm{a},\Gamma,F) =  P\left[{S\over\sigma} \le z\right],\ z \in \R.
$$
Sufficient conditions for $G$ to be close to $\Phi$ are developed.

There has been some recent work on the Central Limit Theorem for sums of linear processes with values in a Hilbert space \cite{N00}.	Other recent work on the Central Limit Theorem for linear random fields has emphasized the Beveridge Nelson decomposition, developed by Phillips and Solo \cite{PS92} for processes.  This approach leads to functional versions of the CLT, \cite{KKK08}, \cite{P10}, and their references.  Our approach follows that of Ibragimov \cite{I63}.  For the case in which $\Gamma$ is a rectangle, it requires no additional conditions on the coefficients, but does require the innovations to be independent and does not deliver a functional version.

\section{Generalities} 

Let
\be\label{eq:rho}
	\rho = \rho(\bm{a},\Gamma) = \max_{r,s\in\Z} {|b_{r,s}|\over\sigma}.
\ee
Interest in $\rho$ stems from the following:

\begin{prop}\label{prop:ksbnd} 
	Let ${\cal H}$ denote a class of distribution functions for which
\be\label{eq:aah}
	\int_{\R} xH\{dx\} = 0,\quad \int_{\R} x^2H\{dx\} = 1, {\rm for\ all}\ H \in {\cal H},
\ee
$$
	\lim_{c\to\infty}\sup_{H\in{\cal H}} \int_{|x|> c} x^2H\{dx\} = 0;\eqno(\ref{eq:aah})
$$
Then $\forall\ \epsilon > 0$, $\exists\ \delta = \delta_{\epsilon,\cal H}$, depending only on $\epsilon$ and ${\cal H}$ for which for which
\be
d(G, \Phi) := \sup_z |G(z)-\Phi(z)| \le \epsilon
\label{eq:ks}
\ee 
for all all $F \in {\cal H}$ for all arrays $\bm{a}$ and regions $\Gamma\subset \Z^2$ for which $\rho \le \delta$.
\end{prop}


	{\it Proof}. Let $\hat{\ }$  denote Fourier transform (characteristic function), so that $\hat{F}(t) = \int_{\R} e^{\imath tx}F\{dx\}$, and 
$$
	L(\eta) = {1\over\sigma^2}\sum_{r,s\in\Z} \int_{|b_{r,s}x|>\eta\sigma} |b_{r,s}x|^2F\{dx\}
$$
for $\eta > 0$.  Then, for any $\eta > 0$, $|\hat{G}(t)-\hat{\Phi}(t)| \le \eta|t|^3/6 + t^2L(\eta)$ for all $t \in \R$ from the proof the the Central Limit Theorem for independent summand (\cite{BL66}, pp. 359 - 361), and 
\begin{align}\label{eq:ksbnd}
	\sup_z |G(z)-\Phi(z)| &\le {1\over\pi} \int_{-T}^T \left|{\hat{G}(t)-\hat{\Phi}(t)\over t}\right|dt + {24 \over \pi \sqrt{2 \pi } T} \\
		&\le  {\eta T^3\over 18} + {1\over 2}T^2L(\eta) + {24 \over T} \nonumber
\end{align}
for any $T > 0$ by the smoothing inequality (\cite{F66}, pp. 510 - 512).  Given $\epsilon$, let $T_{\epsilon} = {96 \over \epsilon} $ and $\eta_{\epsilon} = 4\epsilon/T_{\epsilon}^3$.   Then the left side of (\ref{eq:ksbnd}) is at most ${1\over 2}T_{\epsilon}^2L(\eta_{\epsilon}) + {1\over 2}\epsilon$. Next, let $J(c) = \sup_{H\in{\cal H}} \int_{|x|>c} x^2H\{dx\}$ for $c > 0$, so that $J(c) \to 0$ as $c \to \infty$ by (\ref{eq:aah}); and let $J^{\#}(z) = \inf\{c > 0: J(c) \le z\}$ for $z > 0$.  Then 
$$
	L(\eta) \le \sum_{r \in \Z} \sum_{s\in \Z} b_{r,s}^2\int_{|x|>\eta/\rho} x^2F\{dx\} \le J\left({\eta\over\rho}\right),
$$ 
and $\delta = \eta_{\epsilon}/J^{\#}(T_{\epsilon}^{-2}\epsilon)$ has the desired properties. \qed

	Next, let
$$
		\Vert\bm{a}\Vert_p = \left[\sum_{r,s\in\Z} |a_{r,s}|^p\right]^{1\over p} \le \infty
$$
for $1 \le p \le 2$.  Thus, $\Vert\bm{a}\Vert_2$ is assumed to be finite and $\Vert\bm{a}\Vert_p$ may be finite for some value of $p < 2$.  In terms of $\Vert\bm{a}\Vert_p$ there is a simple bound on $\rho$,
\be\label{eq:crude}
	\rho \le {\Vert\bm{a}\Vert_p\time\#\Gamma^{1\over q}\over\sigma},
\ee
where $q$ denotes the conjugate, $1/p + 1/q = 1$ and $\#\Gamma$ denotes the cardinality of $\Gamma$.  In particular, $\rho \le \Vert\bm{a}\Vert_1/\sigma$.
This leads to:

\begin{cor} Let ${\cal H}$ be as in Proposition 1.  If $\Vert\bm{a}\Vert_1 < \infty$, then $\forall\ \epsilon > 0,\ \exists\ \kappa = \kappa(\epsilon,{\cal H}) > 0$ for which (\ref{eq:ks}) holds whenever $\sigma \ge \kappa\Vert\bm{a}\Vert_1$ and $F \in {\cal H}$.
\end{cor}

	{\it Proof}. For $\delta$ as per in { \it Proposition 1}, pick $\kappa$ such that $ \kappa \delta \leq 1 $. The result is then an easy consequence of the proposition  and the fact that   $\rho \le \Vert\bm{a}\Vert_1/\sigma$. \qed

\section{Rectangles}

	To bound $\rho$, suppose that the maximum occurs when $r = r_0$ and $s = s_0$, say
$$
		|b_{r_0, s_0}| = \max_{r,s} |b_{r,s}|,
$$
and let $\Delta b_{u,v} = b_{u,v}-b_{u,v-1}-b_{u-1,v}+b_{u-1,v-1}$ for $(u,v) \in \Z^2$.  Then
\be\label{eq:Dlta1}
	b_{r_0+r,s_0+s} - b_{r_0,s_0+s} -b_{r_0+r,s_0} + b_{r_0,s_0} = \sum_{r=r_0+1}^{r_0+r}\sum_{v=s_0+1}^{s_0+s} \Delta b_{u,v}
\ee
for $r,s \ge 1$.  Let
\be\label{eq:que}
	Q_{m,n}  = \sum_{r=1}^m \sum_{s=1}^{n} \sum_{u=r_0+1}^{r_0+r}\sum_{v=s_0+1}^{s_0+s} |\Delta b_{u,v}| = \sum_{r=r_0+1}^{r_0+m} \sum_{s=s_0+1}^{s_0+n} (r-r_0)(s-s_0)|\Delta{b}_{r,s}|
\ee
for $m,n \ge 1$.  Then $|b_{r_0,s_0}| \le |b_{r_0,s_0+s}| + |b_{r_0+r,s_0}| + |b_{r_0,s_0}| = \sum_{r=r_0+1}^{r_0+r}\sum_{v=s_0+1}^{s_0+s} |\Delta b_{u,v}|$ for all $r,s \ge 1$ and, therefore,
$$
	mn|b_{r_0,s_0}| \le \sum_{r=1}^m \sum_{s=1}^{n} (|b_{r_0+r,s_0+s}| + |b_{r_0,s_0+s}| + |b_{r_0+r,s_0}|) + Q_{m,n}
$$
for all $m,n \ge 1$.  Here
$$
	\sum_{r=1}^m \sum_{s=1}^{n} |b_{r_0+r,s_0+s}| \le \sqrt{mn}\sqrt{\sum_{r=1}^m \sum_{s=1}^{n} b_{r_0+r,s_0+s}^2} \le \sqrt{mn}\sigma,
$$
and similarly, $\sum_{r=1}^m \sum_{s=1}^{n} |b_{r_0,s_0+s}| \le m\sqrt{n}\sigma$ and $\sum_{r=1}^m \sum_{s=1}^{n} |b_{r_0+r,s_0}| \le \sqrt{m}n\sigma$.  So, $mn|b_{r_0,s_0}| \le \sqrt{mn}\sigma + m\sqrt{n}\sigma + \sqrt{m}n\sigma + Q_{m,n}$.
That is,
\be\label{eq:bscbnd}
	\rho = {|b_{r_0,s_0}|\over\sigma} \le ({2\over\sqrt{m}} + {2\over\sqrt{n}}) + {Q_{m,n}\over mn\sigma}
\ee
for any $m,n \ge 1$.  The first two terms may be made small by taking $m$ and $n$ large.  Thus, the issue is $Q_{m,n}$.  Suppose now that $\Gamma$ can be written as the union of $\ell$ non-empty pairwise mutually exclusive rectangles,
\be\label{eq:unnrctngl}
	\Gamma = \bigcup_{i=1}^{\ell} \{(j,k): \underline{M}_i \le j \le \overline{M}_i,\ \underline{N}_i \le k \le \overline{N}_i\}.
\ee

\begin{prop}
	If  $\Gamma$ is of the form (\ref{eq:unnrctngl}), then 
where $M,N \ge 1$, then 
\be\label{eq:rhorctngle}
	\rho \le 12\left({\sqrt{\ell} \Vert\bm{a}\Vert_2\over\sigma}\right)^{1\over 5} + {8\sqrt{\ell} \Vert\bm{a}\Vert_2\over\sigma}.
\ee
\end{prop}

	{\it Proof}.  In this case $b_{r,s} = \sum_{i=1}^{\ell} b_{r,s}^{(i)}$, where $b_{r,s}^{(i)} = \sum_{j=\underline{M}_i}^{\overline{M}_i} \sum_{k=\underline{N}_i}^{\overline{N}_i} a_{j+r,k+s},\ \Delta{b}_{r,s} = \sum_{i=1}^{\ell} \Delta b_{r,s}^{(i)}$, and $\Delta b_{r,s}^{(i)} = a_{r+\overline{M}_i,s+\overline{N}_i} - a_{r+\underline{M}_i,s+\overline{N}_i} - a_{r+\overline{M}_i,s+\underline{N}_i} + a_{r+\underline{M}_i,s+\underline{N}_i}$.  So, 
$$
	\eqalign{Q_{m,n} &\le  mn\sum_{r=1}^m \sum_{s=1}^{n} \sum_{i=1}^{\ell}  (|a_{r+\overline{M}_i,s+\overline{N}_i} - a_{r+\underline{M}_i,s+\overline{N}_i} - a_{r+\overline{M}_i,s+\underline{N}_i} + a_{r+\underline{M}_i,s+\underline{N}_i}|)\cr
		&\le 4(mn)^2\sqrt{\ell}\Vert\bm{a}\Vert_2\cr},
$$
by Schwartz' Inequality, and 
$$
	\rho \le ({2\over\sqrt{m}} + {2\over\sqrt{n}}) + {4mn\sqrt{\ell}\Vert\bm{a}\Vert_2\over\sigma}
$$
for any $m,n \ge 1$.  Letting $m = n = \lceil ({\sigma/\sqrt{\ell}\Vert\bm{a}\Vert_2})^{4\over 5}\rceil$, the least integer that exceeds $(\sigma/\sqrt{\ell}\Vert\bm{a}\Vert_2)^{4\over 5}$ then leads to (\ref{eq:rhorctngle}). \qed

	When specialize to (intersections of) rectangles (with $\Z^2$), the proposition provides a complete analogue of Ibragimov's theorem \cite{I63} with lots of uniformity.

\begin{cor}.  Let ${\cal H}$ be as in Proposition \ref{prop:ksbnd} and let ${\cal R}_{\kappa}$  be the collections of pairs $(\bm{a},\Gamma)$ for which $\Vert\bm{a}\Vert_2 > 0$, $\Gamma$ is the a rectangle, and $\sigma(\bm{a},\Gamma) \ge \kappa\Vert\bm{a}\Vert_2$.  Then, as $\kappa \to \infty$, the distributions of $S/\sigma$ converge to the $\Phi$ uniformly with respect to $(\bm{a},\Gamma) \in {\cal R}_{\kappa}$ and $F \in {\cal H}$
\end{cor}

\ed